\documentclass{article}
\usepackage{amsmath,amsthm,amsfonts,eucal}

\numberwithin{equation}{section}

\def\cf{{\mathcal F}}

\def\cai{{\mathcal I}}

\def\ck{{\mathcal K}}

\def\cam{{\mathcal M}}

\def\br{{\mathbb R}}

\def\d{\delta}        
\def\eps{\varepsilon}

\def\l{\lambda}       
\def\m{\mu}
\def\n{\nu}

\def\t{\tau}

\def\imply{\Rightarrow}

\def\eq{\Leftrightarrow}

\def\itm#1{\item[($#1$)]}

\def\subc{\underline{\delta}}
\def\supc{\overline{\delta}}

\def\sinc{S^{\uparrow}}
\def\sdec{S^{\downarrow}}

\newtheorem{Thm}{Theorem}[section]
\newtheorem{Cor}[Thm]{Corollary}
\newtheorem{Prop}[Thm]{Proposition}
\newtheorem{Lemma}[Thm]{Lemma}
\theoremstyle{definition}

\theoremstyle{remark}
\newtheorem{rem}[Thm]{Remark} 
 
\begin{document}
\title{\huge  On the domain of singular traces.}
\author{Daniele Guido$^{*}$, Tommaso Isola
\thanks{Supported in part by GNAMPA}
}
\date{March 4, 2002}
\markboth{On the domain of singular traces.}
{On the domain of singular traces.}
\maketitle
\bigskip\bigskip\noindent
Dipartimento di Matematica, Universit\`a di Roma ``Tor
Vergata'', I--00133 Roma, Italy. E-mail: {\tt guido@mat.uniroma2.it, 
isola@mat.uniroma2.it}

\section{Introduction.}

After the introduction of spectral triples in Alain Connes'
Noncommutative Geometry, singular traces on $B(H)$ became a quite
popular tool in operator algebras.  With the aim of classifying them
some papers have been written
(\cite{Dixmier},\cite{Varga},\cite{AGPS}, see also \cite{DFWW} for
nonpositive traces), addressing in particular the question of singular
traceability, namely when, for a given operator, there exists a
singular trace taking a finite non-zero value on it.  Such a problem
has been completely solved in \cite{AGPS}, in terms of a condition on
the asymptotics of the eigenvalue sequence of the given operator, thus
selecting a class of operators whose eigenvalue asymptotics are in
some sense close to the sequence $1/n$.

The problem we address in this note is slightly different, requesting 
now that the given operator just belongs to the domain of a singular 
trace, the value zero being allowed.  Also, a dual question is 
addressed, namely for which operators there exists a singular trace 
assuming the value $\infty$ on it.

A na\"\i ve answer would say that too slow asymptotics cannot be
resummed by any (singular) trace, hence no singular trace may vanish
on them, and conversely that too fast asymptotics are in the kernel of
any (singular) trace, hence no singular trace can be infinite on them. 
Even though this statement will be true for singular traces generated
by `regular' asymptotics (see below), it is false in general.  Indeed
we show that every compact operator is in the kernel (hence in the
domain) of some singular trace, and any infinite-rank operator is not
in the domain (hence not in the kernel) of some singular trace.

The mentioned results will be proved for general semifinite factors,
in the spirit of \cite{GuIs1}.  The proofs are based on the new
characterization of singular traceability given in \cite{GuIs9}, which
we extend here to the continuous case, and also on the fact that
ideals can be described mainly in terms of order relations.

\section{Ideals and singular traces in a semifinite factor}

Here $(\cam,tr)$ is any $\sigma$-finite, semifinite factor with a 
semifinite normal faithful trace, although only in the infinite case 
the following discussion is non trivial.

Let us recall that a function $\m$ is associated with any operator in
$\cam$, via nonincreasing rearrangement (cf.  \cite{FK}): set
$\m_A(t):=\inf\{s\geq0: \l_A(s)\leq t\}$, $t\ge0$, where
$\l_A(t):=tr(e_{|A|}(t,\infty))$, and $|A|=\int_0^\infty t\
de_{|A|}(t)$ is the spectral decomposition of $|A|$.  Recall that
$\m_{A}$ is non-increasing and right continuous, $A$ is
(Breuer-)compact if $\m_{A}$ is infinitesimal, and finite-rank if
$\m_{A}$ is eventually zero.  We also set $g_{A}(t) =
-\log\m_{A}(e^{t})$.

\subsection{Singular traceability}

The singular traceability for a compact operator in $B(H)$, namely the
existence of a singular trace on $B(H)$ which is non-trivial on (the
ideal generated by) $T$, has been completely characterized in
\cite{AGPS}, and this result has been extended to semifinite factors
in \cite{GuIs1}.  Now we will use the Matuszewska indices to give an
alternate description of singular traceability.  Indeed Theorems
\ref{singtrac2} and \ref{singtrac3} have been proved in \cite{GuIs9}
for the case of $B(H)$, but all the arguments extend to the general
factor case.  We give the proofs for the sake of completeness.

Given an operator $A\in\cam$, we define its integral eigenvalue
function $S_{A}$ as
\begin{equation*}
    S_{A}(x) = 
    \begin{cases}
	\sinc_{A}(x):=\int_0^x \m_{A}(y)dy
	&\m_{A}\not\in L^1[0,\infty)\cr
	\sdec_{A}(x):=\int_x^\infty\m_{A}(y)dy
	&\m_{A}\in L^1[0,\infty).
    \end{cases}
\end{equation*}

Recall that a singular trace on $\cam$ is a tracial weight vanishing on 
finite rank operators.

\begin{Thm}\label{singtrac1}
    {\rm \cite{GuIs1}} An operator $T\in\cam$ is singularly traceable if and
    only if $1$ is a limit point, when $x\to\infty$, of the function $
    \frac{S_{A}(\l x)}{S_{A}(x)}, $ for some $\l>1$.  If it is true
    for one $\l$, it is indeed true for any $\l>1$.
\end{Thm}

The singular traceability condition can be reformulated as follows.
 
\begin{Prop} \label{singtrac2}
    $A$ is singularly traceable if and only if
    $\liminf\frac{x\m_{A}(x)}{S_{A}(x)}=0$.
\end{Prop}

\begin{proof}
    Assume first $A$ is not trace class, i.e. $S_{A}(x)=\sinc_{A}(x)$. 
    Then the thesis follows by Theorem \ref{singtrac1} and the
    following inequalities:
    \begin{equation*}
	0<\frac{\sinc_{A}(2x)}{\sinc_{A}(x)}-1
	\leq\frac{x\m_{A}(x)}{\sinc_{A}(x)}
	\leq 2\left(1-\frac{\sinc_{A}(x/2)}{\sinc_{A}(x)}\right).
    \end{equation*}
    When $A$ is trace class, i.e. $S_{A}(x)=\sdec_{A}(x)$, we have,
    analogously,
    \begin{equation*}
	0<1-\frac{\sdec_{A}(2x)}{\sdec_{A}(x)}
	\leq\frac{x\m_{A}(x)}{\sdec_{A}(x)}
	\leq 2\left(\frac{\sdec_{A}(x/2)}{\sdec_{A}(x)}-1\right)
    \end{equation*}
    and the thesis follows.
\end{proof}

We now define the Matuszewska indices $\subc(A)$, $\supc(A)$ for a
compact operator $A$ as the indices for the corresponding eigenvalue
function, cf. \cite{BGT}.  As a consequence,
\begin{align}
	\frac{1}{\subc(A)} 
	&=
	\lim_{h\to\infty}\limsup_{t\to\infty}
	\frac{g_{A}(t+h)-g_{A}(t)}{h}
	=
	\inf_{h>0}\limsup_{t\to\infty}
	\frac{g_{A}(t+h)-g_{A}(t)}{h}
	\label{Matu2}
	\\
	\frac1{\supc(A)} 
	&=\lim_{h\to\infty}\liminf_{t\to\infty}
	\frac{g_{A}(t+h)-g_{A}(t)}{h}=\sup_{h>0}\liminf_{t\to\infty}
	\frac{g_{A}(t+h)-g_{A}(t)}{h}.\label{Matu1}
\end{align}

For the existence of the limits and the equalities in the definition 
above, see e.g. \cite{BGT}. The following Lemma holds.

\begin{Lemma}\label{trace}
    \item{$(i)$} If $\subc(A)>1$ then $g_{A}(t) < c + (1-\eps)t$, for
    suitable $c,\eps>0$.  In particular $A$ is not trace class. 
    \item{$(ii)$} If $\supc(A)<1$ then $g_{A}(t) > -c + (1+\eps)t$, for
    suitable $c,\eps>0$.  In particular $A$ is trace class. 
    \item{$(iii)$} If $\subc(A)=\supc(A)=1$ then $-c_{1}+(1-\eps)t\leq
    g_{A}(t) \leq c_{2}+(1+\eps)t$, for suitable $c_{1},c_{2},\eps>0$.
\end{Lemma}
\begin{proof}
    Equalities (\ref{Matu2}), (\ref{Matu1}) imply
    $$
	\subc(A)^{-1}\geq\limsup_{t\to\infty}
	\frac{g_{A}(t)}{t}\geq\liminf_{t\to\infty}
	\frac{g_{A}(t)}{t}\geq\supc(A)^{-1},
    $$
    from which the thesis follows.
\end{proof}

\begin{Thm}\label{singtrac3}
    Let $A\in\cam$ be compact.  Then $A$ is singularly traceable if
    and only if $\subc(A)\leq1\leq\supc(A)$.
\end{Thm}

\begin{proof}
     Let $\subc(A)>1$.  By (\ref{Matu2}), this is equivalent to the
     existence of $h>0$ such that $ \limsup_{t\to\infty}
     \frac{g_{A}(t+h)-g_{A}(t)}{h} <1$, or, equivalently, to the
     existence of $\l>1$ for which $
     \l\liminf_{t\to\infty}\frac{\m_{A}(\l t)}{\m_{A}(t)}>1.  $ Now
     observe that, by Lemma \ref{trace}, $A$ is not trace class. 
     Therefore
     $$
     \frac{S_{A}(\l x)}{S_{A}(x)} 
     =\frac{\l\int_{0}^{x}\m_{A}(\l t) dt}{\int_{0}^{x}\m_{A}(t) dt} 
     =\l\frac{\int_{0}^{x}\left(\frac{\m_{A}(\l t)}{\m_{A}(t)}\right) 
     \m_{A}(t) dt}{\int_{0}^{x}\m_{A}(t) dt},
     $$
     hence
     $$
     \liminf_{x\to\infty}\frac{S_{A}(\l x)}{S_{A}(x)} \geq 
     \l\liminf_{x\to\infty}\left(\frac{\m_{A}(\l x)}{\m_{A}(x)}\right)>1,
     $$
     which implies that $A$ is not singularly traceable by Theorem
     \ref{singtrac1}.  \\
     The proof when $\supc(A)<1$ is analogous.
     \\
     Assume now that $A$ is not singularly traceable, namely, by
     Proposition \ref{singtrac2},
     $
     \inf_{t}\frac{t\m_{A}(t)}{S_{A}(t)}= k>0$.  If $A$ is not
     trace-class, i.e. $S_{A}=\sinc_{A}$, since $\m_{A}$ is the
     derivative of $\sinc_{A}$, the hypothesis means that
     $$
     \frac{d}{dt} \log\sinc_{A}(t)\geq\frac{k}{t}, \forall t.
     $$
     Integrating on the interval $[x,\l x]$ one gets $ \sinc_{A}(\l
     x)\geq \l^{k}\sinc_{A}(x).  $ Since
     $\frac{x\m_{A}(x)}{\sinc_{A}(x)}\leq 1$, one obtains
     \begin{equation*}
	 \frac{\m_{A}(x)}{\m_{A}(\l
	 x)}\leq\frac{\frac{\sinc_{A}(x)}{x}}{\frac{k\sinc_{A}(\l
	 x)}{\l x}} \leq \frac{\l^{1-k}}{k}.
     \end{equation*}
     As a consequence,
     \begin{align*}
	 \subc(A)^{-1}&=\lim_{\l\to\infty}\frac{1}{\log 
	 \l}\limsup_{x\to\infty}\log\frac{\m_{A}(x)}{\m_{A}(\l x)}\cr 
	 &\leq\lim_{\l\to\infty}\frac{-\log 
	 k+(1-k)\log\l}{\log\l}=1-k<1
     \end{align*}
     \\
     If $A$ is trace-class, i.e. $S_{A}=\sdec_{A}$, since $-\m_{A}$ is
     the derivative of $\sdec_{A}$, we may prove, in analogy with the
     previous argument, that $ \sdec_{A}(\l x)\leq
     \l^{-k}\sdec_{A}(x).  $ Since
     $\sdec_{A}(t)\geq\int_{t}^{2t}\m_{A}(s)ds\geq t\m_{A}(2t)$, which
     implies $\frac{x\m_{A}(2x)}{\sdec_{A}(x)}\leq 1$, one obtains
     \begin{equation*}
	 \frac{\m_{A}(2\l x)}{\m_{A}(x)}\leq\frac{\frac{\sdec_{A}(\l x)}{\l 
	 x}}{\frac{k\sdec_{A}(x)}{x}} \leq \frac{\l^{-1-k}}{k},
     \end{equation*}
     namely
     \begin{equation*}
	 \frac{\m_{A}(x)}{\m_{A}(\l x)} \geq k
	 \left(\frac{\l}{2}\right)^{1+k},
     \end{equation*}
     As a consequence,
     \begin{align*}
	 \supc(A)^{-1}&=\lim_{\l\to\infty}\frac{1}{\log 
	 \l}\liminf_{x\to\infty}\log\frac{\m_{A}(x)}{\m_{A}(\l x)}\cr 
	 &\geq\lim_{\l\to\infty}\frac{\log 
	 k+(1+k)\log2+(1+k)\log\l}{\log\l}=1+k>1.
     \end{align*}
 \end{proof}
 
\begin{rem}We say that $A$ is {\it regular} if $\subc(A)=\supc(A)=:\d(A)$.  As a
consequence, for a regular $A$, singular traceability is equivalent to
$\d(A)=1$.\end{rem}

\subsection{Ideals}

Let us introduce the set $M$ of non-increasing infinitesimal right
continuous functions defined on the interval $[0,\infty)$, and the set
$G$ of $(-\infty,+\infty]$-valued functions on $\br$ which are
non-decreasing, right continuous, bounded from below and unbounded
from above.  Clearly the map
$$
\m\to g(t)=-\log\m(e^{t})
$$
gives an order-reversing one-to-one correspondence from $M$ to $G$. 

Consider the action $\l\to D_{\l}f$ of the multiplicative
group $\br_+$ on $M$ given by:
 $$
D_{\l}f(t)=\l f(\l t)\ ,\qquad\l,t\in\br_+
 $$
We say that a face $F$ in $M$, namely a hereditary subcone of $M$ is
{\it dilation invariant} if $f\in F$ $\imply$ $D_{\l}f\in F$,
$\l\in\br_+$.

We proved in \cite{GuIs1} that proper two-sided ideals in $\cam$ are in one-to-one 
correspondence with dilation invariant faces in $M$, namely for any 
ideal $\cai$ there exists a face $F$ such that
$$
A\in\cai\eq\m_{A}\in F,
$$
and conversely if $F$ is a dilation invariant face then 
$\{A\in\cam:\m_{A}\in F\}$ is a two sided ideal in $\cam$.

We note here that the additivity property of $F$ does not really 
matter, indeed since $\m\vee\n\leq\m+\n\leq2(\m\vee\n)$, homogeneity 
and the closure under $\vee$ are equivalent to linearity.  As a 
consequence ideals of $\cam$ can be described by subsets of $G$, namely 
the following corollary holds.
 
\begin{Cor}\label{F-ideals}
	There is a one to one correspondence between ideals in $\cam$ and 
	subsets $H$ of $G$ such that
	\item $f,g\in H\imply f\wedge g\in H$
	\item $f\in H$ and $g\in G$ with $g\geq f\imply g\in H$
	\item $f\in H \imply f+a\in H$, $\forall a\in\br$
	\item $f\in H \imply f^{a}\in H$, $\forall a\in\br$, where 
	$f^{a}(t)=f(t-a)$. \\
	In particular, given such an $H$, the set 
	$\{A\in\cam:g_{A}\in H\}$ is a two sided ideal in $\cam$.
\end{Cor}

According to the previous correspondence, the ideal $\ck(\cam)$ of 
compact operators corresponds to the whole set $G$, while the 
ideal $\cf(\cam)$ of finite rank operators corresponds to eventually 
infinite elements of $G$. The principal ideal generated by an 
operator $B$ corresponds to the set $H(B)$ defined as follows:
\begin{equation}\label{princIde}
f\in H(B)\Leftrightarrow \exists a,b\in\br: f\geq b+g^{a}_{B}.
\end{equation}

We now introduce the notion of kernel of an ideal.  If $\cai$ is an
ideal in $\cam$, we define its kernel $\cai_{0}$ as the set of
$A\in\cam$ for which there exists $T\in\cai_{+}$ such that
$$\forall\eps>0\ \exists F_{\eps}\in\cf(\cam):
|A|<\eps T+F_{\eps}.$$

The following statement is an immediate consequence of the definition
of $\cai_{0}$.

\begin{Prop} 
	Let $\cai$ be an ideal in $\cam$.  Then $\cai_{0}$ is an ideal 
	contained in $\cai$, and any singular trace defined on $\cai$ 
	vanishes on $\cai_{0}$.
\end{Prop}

Clearly, if the face $F$ corresponds to $\cai$, the face $F_{0}$
corresponding to $\cai_{0}$ is defined as the set of $\m\in M$ for
which there exists $\n\in F$ such that
$$\forall\eps>0\ \exists x_{0}>0: \m(x)<\eps\n(x),\ x>x_{0},$$
hence the  subset $H_{0}$ of $G$, corresponding to $F_{0}$, consists of the
elements $g\in G$ for which there exists $h\in H$ such that 
\begin{equation*}
\forall c>0\ \exists t_{0}\in\br : g(t)>c+h(t),\ t>t_{0}.
\end{equation*}

If $H=H(B)$ corresponds to the ideal generated by $B$, then $g\in 
H_{0}(B)$ if
\begin{equation}\label{ker(B)}
	\exists a\in\br:\forall c>0 \exists t_{0}\in\br: 
	g(t)>c+g_{B}^{a}(t),\ t>t_{0}.
\end{equation}

Let us prove the following.

\begin{Lemma}\label{regideal}
	Let $B$ be regular.  Then a compact operator $A$ is in the ideal
	$\cai(B)$ generated by $B$ if and only if there exists
	$T\in\cai(B)$ regular such that $g_{A}\geq g_{T}$.
\end{Lemma}

\begin{proof}
    By equation (\ref{princIde}), $A\in\cai(B)$ if there exist 
    $a,b\in\br$ such that $g_{A}\geq b+g_{B}^{a}$. Clearly
    $\subc(g_{B})=\subc(b+g_{B}^{a})$, and 
    $\supc(g_{B})=\supc(b+g_{B}^{a})$, therefore, choosing $T$ such 
    that $g_{T}=b+g_{B}^{a}$, we obtain that $T$ is regular.
\end{proof}

\section{Domains and kernels of singular traces}

\begin{Thm}\label{Thm1}
	\itm{i} Let $A$ be a compact operator.  Then there is a singular trace 
	vanishing on $A$.  
	\itm{ii} Let $A$ be an infinite rank operator.  Then there is a singular 
	trace which is infinite on $A$.
\end{Thm}
\begin{proof}
    In this proof, opposed to the notation of Corollary 
    \ref{F-ideals}, the upper indices of $g_{A}$ are true exponents.
    \\
    $(i)$ The statement is proved if we show that there exists a
    singularly traceable operator $B$ such that $A\in\cai_{0}(B)$. 
    According to the previous discussion this amounts to find an
    element $B$ such that $g_{A}$ satisfies condition (\ref{ker(B)}). 
    \\
    Choose inductively an increasing sequence $t_{n}$ such that
    $t_{n+1}-t_{n}>n$ and $g^{1/2}_{A}(t_{n+1})-g^{1/2}_{A}(t_{n})>n$,
    and set $g(t)=g^{1/2}_{A}(t_{n})$ when $t\in[t_{n},t_{n+1})$. 
    Obviously $g\in G$.  Then choose $B$ such that $g_{B}=g$.  Since
    $g\leq g^{1/2}_{A}$ and $\forall c>0$ $\exists t_{0}\in\br :
    g_{A}(t)>c+g^{1/2}_{A}(t)$, $t>t_{0}$, then a fortiori $g_{A}$
    satisfies (\ref{ker(B)}).  Moreover one easily gets $\subc(g)=0$
    and $\supc(g)=\infty$, which implies that $B$ is singularly
    traceable.  \\
    $(ii)$ The statement is proved if we show that there exists a
    singularly traceable operator $B$ such that $A\not\in\cai(B)$.  In
    analogy with the previous proof, this amounts to find a function
    $g\in G$ such that
    \begin{equation}\label{infinite}
	\forall c>0\ \exists t_{0}\in\br: g_{A}(t)<c+g(t),\ t>t_{0},
    \end{equation}
    and then choosing $B$ such that $g_{B}=g$.
    \\
    Choose inductively an increasing sequence $t_{n}$ such that
    $t_{n+1}-t_{n}>n$ and $g^{2}_{A}(t_{n+1})-g^{2}_{A}(t_{n})>n$, and
    set $g(t)=g^{2}_{A}(t_{n+1})$ when $t\in[t_{n},t_{n+1})$. 
    Obviously $g\in G$.  Since $g^{2}_{A}$ clearly satisfies
    (\ref{infinite}) and $g\geq g^{2}_{A}$, a fortiori $g$ satisfies
    (\ref{infinite}).  Moreover one easily gets $\subc(g)=0$ and
    $\supc(g)=\infty$, which implies that $B$ is singularly traceable.
\end{proof}

\begin{Thm}\label{Thm2}
	Let $A$ be not singularly traceable, and $B$ satisfy
	$\d(B)=1$.  
	\itm{i} If $A$ is not trace class, then $A$ does not belong to
	the ideal $\cai(B)$ generated by $B$, namely any singular
	trace on $\cai(B)$ is infinite on $A$. 
	\itm{ii} If $A$ is trace class, then $A$ belongs to the kernel
	$\cai_{0}(B)$ of $\cai(B)$, namely any singular trace on
	$\cai(B)$ is zero on $A$.
\end{Thm}
\begin{proof}
    $(i)$ By Theorem \ref{singtrac3}, $\subc(A)>1$, hence, by Lemma
    \ref{trace} there exist $\eps>0$ and $t_{0}>0$ such that
    $g_{A}(t)\leq (1-\eps)t$ when $t>t_{0}$.  Now, for any regular
    $T\in\cai(B)$, we have, by Lemma \ref{trace}, $g_{T}(t)\geq c 
    (1-\eps/2)t$. 
    This implies that eventually $g_{T}(t)\geq g_{A}(t)$.  The thesis
    then follows from Lemma \ref{regideal}.  \\
    $(ii)$ By Theorem \ref{singtrac3}, $\supc(A)<1$, hence, by Lemma
    \ref{trace}, there exist $\eps>0$ and $t_{0}>0$ such that
    $g_{A}(t)\geq (1+\eps)t$ when $t>t_{0}$.  Now, we have,
    by Lemma \ref{trace}, 
    $g_{B}(t)\leq c (1+\eps/2)t$.  This clearly implies $g_{A}$
    satisfies property (\ref{ker(B)}) with $a=0$, hence the thesis.
\end{proof}

 

\begin{thebibliography}{99}

 \bibitem{AGPS} S. Albeverio, D.~Guido, A.~Ponosov, S.~Scarlatti.  
 {\it Singular traces and compact operators}.  J. Funct.  Anal., {\bf 
 137} (1996), 281--302.
   
 \bibitem{BGT} N.H. Bingham, C.M. Goldie, J.L. Teugels. {\it Regular 
 variation}. Cambridge University Press, Cambridge, 1987.

 \bibitem{Dixmier} J. Dixmier.  {\it Existence de traces non 
 normales}.  C.R. Acad.  Sci.  Paris, {\bf 262} (1966), 1107--1108.
 
 \bibitem{DFWW} K. Dykema, T. Figiel, G. Weiss, M. Wodzicki.  {\it 
 Commutator Structure of Operator Ideals}, preprint MSRI  n. 2001-013.

 \bibitem{FK} T. Fack, H. Kosaki.  {\it Generalized s-numbers of 
 $\t$-measurable operators}.  Pacific J. Math., {\bf 123} (1986), 269.

 \bibitem{GuIs1} D. Guido, T. Isola.  {\it Singular traces for 
 semifinite von~Neumann algebras}.  J. Funct.  Anal.,
 {\bf 134} (1995), 451--485.

 \bibitem{GuIs9} D. Guido, T. Isola.  {\it Dimension and singular traces 
 for spectral triples, with applications to fractals}, preprint  math.OA/0202108.

 \bibitem{Varga} J.V. Varga.  {\it Traces on irregular ideals}.  Proc.  
 A.M.S., {\bf107} (1989), 715.

 \end{thebibliography}
 \end{document}